\documentclass[a4paper, 10pt]{article}
\usepackage{geometry}
\usepackage[francais, english]{babel}
\usepackage{latexsym,amsfonts, amssymb}
\usepackage{amsthm}
\usepackage[latin1]{inputenc}
\usepackage{amsmath}
\usepackage{amssymb}
\usepackage{amsfonts}
\usepackage{dsfont}
\usepackage{color}
\usepackage{multicol}
\usepackage{hyperref}
\usepackage{amsthm}
\usepackage{latexsym}
\usepackage{fancyhdr}
\usepackage{epsfig}
\usepackage{mathrsfs}
\usepackage{multirow}

\makeatletter
\renewcommand\theequation{\thesection.\arabic{equation}}
\@addtoreset{equation}{section} \makeatother

\newfont{\bbf}{cmbx12 scaled 1435}

\newtheorem{prop}{Proposition}[section]
\newtheorem{lem}{Lemma}[section]

\renewcommand{\theequation}{\thesection.\arabic{equation}}
\newcommand{\eop}{\hspace*{\fill} \ensuremath{\Box}}

\begin{document}

\newcommand{\I}[1]{\mathds{1}_{{#1}}}

\def\Sum{ \displaystyle \sum }
\def\Frac{\displaystyle \frac}

\def\Cit{\mathbb{C}}
\def\esp{\mathbb{E}}
\def\Var{\hbox{\rm Var}}
\def\Cov{\hbox{\rm Cov}}
\def\Supp{\hbox{\rm Supp}}
\def\Card{\hbox{\rm Card}}
\def\Det{\hbox{\rm Det}}
\def\Tr{\hbox{\rm Tr}}
\def\Dim{\hbox{\rm dim}}
\def\Rank{\hbox{\rm dim}}
\def\Id{\hbox{\rm Id}}
\def\Ker{\hbox{\rm Ker}}
\def\ind{\mathbb{I}}
\def\Nit{\mathbb{N}}
\def\Rit{\mathbb{R}}
\def\Zit{\mathbb{Z}}
\def\prob{\mathbb{P}}
\def\where{\rm where}
\def\with{\rm with}
\def\I{\rm I}
\def\J{\rm J}
\def\Ip{\rm I^{\prime}}
\def\Jp{\rm J^{\prime}}
\def\as{\rm a.s}

\setcounter{page}{1} \setlength{\baselineskip}{.32in}

\begin{center}
 {\Large\bf Auxiliary results for
``Nonparametric kernel estimation of the probability density
function of  regression errors using estimated residuals"}
\end{center}

\begin{center}
By Rawane Samb
\\ Université du Québec à Montréal, Canada
\end{center}

\begin{center}
This version: March 2012
\end{center}

\vskip0.3cm\noindent
Let $(X_{1},Y_{1}),\ldots,(X_{n},Y_{n})$ be a sample  of
independent replicates of the random vector $(X, Y)$, where $Y$ is
the univariate dependent variable and $X$ is the covariate of
dimension $d$. Let $m(\cdot)$ be the conditional expectation of
$Y$ given $X$ and let
 $\varepsilon$ be
the related regression error term, so that the regression model is
\begin{eqnarray*}
Y =m(X)+\varepsilon,
\end{eqnarray*}
 where $\varepsilon$ is assumed to have mean zero and to be
 statistically independent of $X$, and the function $m$
 is smooth but unknown. Our aim is to investigate the nonparametric estimation
of the probability density function of the error term
$\varepsilon$.

\section{Construction of the estimator}

Define
$$
\widehat{\varepsilon}_i
=
 Y_i-\widehat{m}_{in},
\quad
i=1,\ldots,n,
$$
where $\widehat{m}_{in} =\widehat{m}_{in}(X_i)$ is the leave-one
out
 version of  the Nadaraya-Watson (1964) kernel estimator of
 $m(X_i)$,
$$
\widehat{m}_{in}
=
 \frac{\sum_{j=1\atop j\neq i}^nY_j
K_0\left(\frac{X_j-X_i}{b_0}\right)}
{\sum_{j=1\atop j\neq i}^n
K_0\left(\frac{X_j-X_i}{b_0}\right)}.
$$
Here $K_0(\cdot)$ is a kernel function defined on $\Rit^d$ and
$b_0=b_0(n)$ is a bandwidth sequence.

The proposed estimator for the density $f$ of $\varepsilon$ is
$$
\widehat{f}_{n}(e)
 =
\frac{1}{b_1\sum_{i=1}^n \mathds{1}
\left(X_i\in\mathcal{X}_0\right)}
\sum_{i=1}^n \mathds{1}
\left(X_i\in \mathcal{X}_0\right)
K_1\left(
\frac{\widehat{\varepsilon}_i-e}{b_1}
\right),
\quad e\in\Rit,
$$
where $\mathcal{X}_0$  is an inner subset of the support
$\mathcal{X}$ of the covariate $X$,
 $K_1(\cdot)$ is a univariate kernel function and
$b_1=b_1(n)$ is a bandwidth sequence.

\section{Assumptions}

 \vskip 0.3cm \noindent
 {$\bf(A_1)$}
 {\it
 The support $\mathcal{X}$ of $X$  is a subset of
 $\Rit^d$, $\mathcal{X}_0$
 has a nonempty interior and the closure of $\mathcal{X}_0$
 is in the interior of $\mathcal{X}$.
}
\medskip
\\{$\bf(A_2)$}
{\it The p.d.f. $g(\cdot)$ of the i.i.d. covariates $X_i$  is
strictly positive over the closure of $\mathcal{X}_0$, and has
continuous second order partial derivatives  over $\mathcal{X}$.}
\medskip
\\{$\bf(A_3)$}
{\it
 The regression function $m(\cdot)$ has continuous second order partial
derivatives  over  $\mathcal{X}$.
}
\medskip
\\{$\bf(A_4)$}
{\it The i.i.d. centered error regression terms
 $\varepsilon_i$ have finite 6th moments and are independent of the
 covariates $X_i$.
}
\\{$\bf(A_5)$}
{\it
 The  probability  density function $f(\cdot)$ of the $\varepsilon_i$'s
 has  bounded continuous second order
 derivatives over  $\Rit$ and satisfies
 $\sup_{e\in\Rit}|h_p^{(k)}(e)|<\infty$, where  $h_p (e) = e^p f(e)$,
 $p\in[0,2]$ and $k\in\{0,1,2\}$.
  }
\\{$\bf(A_6)$}
{\it The kernel function $K_0$  is  symmetric, continuous over
$\Rit^d$ with support  contained in $[-1/2, 1/2]^d$ and satisfies
$\int\! K_0 (z) dz = 1$.
}
\medskip
\\{$\bf(A_7)$}
{\it
 The kernel function $K_{1}$ is symmetric, has a compact support,  is three times
 continuously differentiable over
 $\Rit$, and satisfies $\int\! K_1 (v) dv = 1$,
 $\int\! K_1^{(\ell)} (v) dv = 0$
 for $\ell=1,2,3$, and
 $\int\! v K_1^{(\ell)}(v) dv=0$ for $\ell=2,3$.
 }
\medskip
\\{$\bf(A_{8})$}
{\it
 The bandwidth $b_0$ decreases to $0$ when $n\rightarrow\infty$ and satisfies,
 for $d^*=\sup\{d+2,2d\}$, $nb_0^{d^*}/\ln n\rightarrow\infty$ and
 $\ln(1/b_0)/\ln(\ln n)\rightarrow\infty$
 when $n\rightarrow\infty$.
}
\\
{$\bf(A_{9})$}
 {\it
The bandwidth $b_1$ decreases to $0$ and satisfies
$n^{(d+8)}b_1^{7(d+4)}\rightarrow\infty$ when
$n\rightarrow\infty$.}

\section{Auxiliary results}

\begin{prop}

Define
$$
\beta_{in}
 =
\frac{\mathds{1} \left(X_i\in\mathcal{X}_0\right)}
{nb_0^d\widehat{g}_{in}} \sum_{j=1, j\neq i}^n
\left(m(X_j)-m(X_i)\right) K_0\left(\frac{X_j-X_i}{b_0}\right),
$$
where
$$
\widehat{g}_{in}
=
\frac{1}{nb_0^d}\sum_{j=1, j\neq i}^n
 K_0\left(\frac{X_j-X_i}{b_0}\right).
$$
Then, under $(A_1)-(A_{9})$, we have, for
all $e\in\Rit$ and $b_0$ and $b_1$ go to $0$, $\sup_{i}|\beta_{in}|=O_{\prob}(b_0^2)$ and
\begin{eqnarray*}
\sum_{i=1}^n \beta_{in}
 K_1^{(1)}
 \left(
 \frac{\varepsilon_i-e}{b_1}
 \right)
 =
 O_{\prob}
 \left(b_0^2\right)
 \left(nb_1^2+(nb_1)^{1/2}\right).
\end{eqnarray*}
\label{Betasum}
\end{prop}

\begin{prop}
Set
$$
\Sigma_{in}
 =
\frac{\mathds{1}\left(X_i\in\mathcal{X}_0\right)}
{nb_0^d\widehat{g}_{in}} \sum_{j=1, j\neq i}^n
 \varepsilon_j
K_0\left(\frac{X_j-X_i}{b_0}\right).
$$
Then, under $(A_1)-(A_{9})$, we have,
for all $e\in\Rit$ and $b_0$ and $b_1$ going to $0$,
\begin{eqnarray*}
\sum_{i=1}^n \Sigma_{in}
 K_1^{(1)}
 \left(
 \frac{\varepsilon_i-e}{b_1}
 \right)
 =
 O_{\prob}
 \left(
 nb_1^4
 +
 \frac{b_1}{b_0^d}
\right)^{1/2}.
\end{eqnarray*}
\label{Sigsum}
\end{prop}

\begin{prop}
Let $\Var_n(\cdot)$ be the
conditional variance given
$X_1,\ldots,X_n$, and set
\begin{eqnarray*}
\zeta_{in}
=
\mathds{1}
\left(X_i \in \mathcal{X}_0\right)
(\widehat{m}_{in}-m(X_i))^2
K_1^{(2)}
\left( \frac{\varepsilon_i-e}{b_1} \right).
\end{eqnarray*}
Then under ${\rm (A_1)-(A_{9})}$, we have, for all $e\in\Rit$, and
$b_0$ and $b_1$ going to $0$,
\begin{eqnarray*}
\Var_n\left(\sum_{i=1}^n\zeta_{in}\right)
 =
O_{\prob}
\left(
nb_1
 +
 n^2b_0^db_1^{7/2}
 \right)
\left(b_0^4+\frac{1}{nb_0^d}\right)^2.
\end{eqnarray*}
\label{VarZn}
\end{prop}

\begin{prop}
Let $R_{in}=\mathds{1}
 \left(X_i\in\mathcal{X}_0\right)
 \left(\widehat{m}_{in}-m(X_i)\right)^3I_{in}$, where
\begin{eqnarray*}
I_{in}
=
\int_{0}^{1}
 (1-t)^2
 K_1^{(3)}
 \left(
 \frac{
 \varepsilon_i-t(\widehat{m}_{in} -m(X_i))-e}{b_1}
 \right)
 dt.
\end{eqnarray*}
Then, under ${\rm (A_1)-(A_{9})}$, we have, for all $e\in\Rit$, and
$b_0$ and $b_1$ going to $0$,
$$
\Var_n\left(\sum_{i=1}^n R_{in}\right)
 =
 O_{\prob}
 \left(n^2b_0^db_1\right)
 \left(b_0^4+\frac{1}{nb_0^d}\right)^3.
$$
\label{VarRn}
\end{prop}

\section{Intermediate results for the propositions}

The proofs of the propositions are based on the following results.

\begin{lem}
Define
$$
\widehat{g}_n(x)
 =
 \frac{1}{nb_0^d}
 \sum_{i=1}^n
 K_0
 \left(
 \frac{X_i-x}{b_0}
 \right),
 \quad
 \overline{g}_n(x)
 =
 \esp\left[\widehat{g}_n(x)\right],
 \quad
 x\in\mathcal{X}_0.
$$
Then under $(A_1)-(A_2)$,  $(A_6)$ and $(A_8)$,
 we have, when $b_0$ goes to $0$,
$$
\sup_{x\in\mathcal{X}_0}
\left|\overline{g}_{n}(x)-g(x)\right|
 =
O\left(b_0^{2}\right),
\quad \sup_{x\in\mathcal{X}_0}
 \left|
\widehat{g}_{n}(x)
 -
 \overline{g}_n(x)
\right|
 = O_{\prob}
 \left(
 b_0^4 +
 \frac{ \ln n}{nb_0^d}
 \right)^{1/2},
$$
 and
$$
\sup_{x\in\mathcal{X}_0}
 \left|
 \frac{1}{\widehat{g}_{n}(x)}
   -
\frac{1}{g(x)} \right|
=
 O_{\prob}
 \left(
 b_0^4
 +
 \frac{\ln n}{nb_0^d}
\right)^{1/2}.
$$
\label{Estig}
\end{lem}

\begin{lem}
Assume that  ${\rm (A_4)}$ and ${\rm (A_6)}$ hold. Then, for any
$1 \leq i \neq j \leq n$,
$$
\left(\widehat{m}_{in}-m(X_i),\varepsilon_i\right)
 \mbox{\it and }
 \left( \widehat{m}_{jn}-m(X_j),\varepsilon_j\right)
$$
are independent given $X_1, \ldots, X_n$,
 provided that $\| X_i - X_j \| \geq C b_0$,
for some  $C>0$.
\label{Indep}
\end{lem}

\begin{lem}

Let $\esp_n[\cdot]$ be the
conditional mean given
$X_1,\ldots,X_n$, and
assume $(A_1)-(A_{9})$. Then,
\begin{eqnarray*}
\sup_{1\leq i\leq n}
 \esp_{n}
  \left[
  \mathds{1}
  \left(X_i \in \mathcal{X}_0\right)
 (\widehat{m}_{in} - m(X_i))^4
 \right]
 &=&
 O_{\prob}
 \left(b_0^4+\frac{1}{nb_0^d}\right)^2,
 \\
 \sup_{1\leq i\leq n}
 \esp_{n}
 \left[
 \mathds{1}
 \left(X_i \in \mathcal{X}_0\right)
 (\widehat{m}_{in} - m(X_i))^6
 \right]
 &=&
 O_{\prob}
 \left(b_0^4+\frac{1}{nb_0^d}\right)^3.
\end{eqnarray*}
\label{BoundEspmchap}
\end{lem}

\begin{lem}

Under  ${\rm (A_5)}$ and ${\rm (A_7)}$ we have, for some $C>0$,
and for any $e\in\Rit$ and $p\in[0,2]$,
\begin{eqnarray}
\left|
 \int
 K_1^{(1)}
 \left( \frac{\epsilon-e}{b_1} \right)^2
 \epsilon^p f(\epsilon)
  d\epsilon
 \right|
 \leq C b_1,
 &&
 \left|
  \int
  K_1^{(1)}
  \left(
 \frac{\epsilon-e}{b_1}
 \right)
 \epsilon^p f(\epsilon)
 d\epsilon
  \right|
 \leq C b_1^2,
\label{MomderK1}
\\
\left|
\int
 K_1^{(2)}
 \left( \frac{\epsilon-e}{b_1} \right)^2
  \epsilon^p f(\epsilon)
  d\epsilon
\right|
 \leq
  C b_1,
  &&
 \left|
\int K_1^{(2)}
 \left(
 \frac{\epsilon-e}{b_1}
 \right)
 \epsilon^p f(\epsilon)
 d\epsilon
\right|
 \leq C b_1^3,
 \label{MomderK2}
\\
\left| \int
 K_1^{(3)}
 \left( \frac{\epsilon-e}{b_1} \right)^2
  \epsilon^p f(\epsilon)
  d\epsilon
\right|
 \leq
  C b_1,
  &&
 \left|
\int K_1^{(3)}
 \left(
 \frac{\epsilon-e}
 {b_1}
 \right)
 \epsilon^p f(\epsilon)
 d\epsilon
\right|
 \leq
 C b_1^3.
 \label{MomderK3}
\end{eqnarray}
\label{MomderK}
\end{lem}

\noindent The proof of all these lemmas is postponed in the appendix.

\section{Proofs of the auxiliary results}

\subsection*{Proof of Proposition \ref{Betasum}}

Assumption $(A_4)$ and Lemma \ref{MomderK}-(\ref{MomderK1}) yield
\begin{eqnarray*}
\left|
\esp_n
 \left[
 \sum_{i=1}^n
 \beta_{in}
 K_1^{(1)}
 \left(
\frac{\varepsilon_i-e}{b_1}
\right)
 \right]
 \right|
 & = &
\left|
 \esp
 \left[
 K_1^{(1)}
  \left(
  \frac{\varepsilon-e}{b_1}
  \right)
  \right]
  \sum_{i=1}^n
  \beta_{in}
   \right|
  \leq C n b_1^2
 \max_{1 \leq i \leq n}
 \left| \beta_{in} \right|,
\\
\Var_n
 \left[
 \sum_{i=1}^n
  \beta_{in}
  K_1^{(1)}
  \left(
\frac{\varepsilon_i-e}{b_1}
 \right)
 \right]
&\leq&
  \sum_{i=1}^n
 \beta_{in}^2
 \esp
 \left[
 K_1^{(1)}
  \left(
   \frac{\varepsilon-e}{b_1}
\right)^2
 \right]
  \leq
 C n b_1
  \max_{1 \leq i \leq n}
  \left|
  \beta_{in}
\right|^2 .
\end{eqnarray*}
Hence the (conditional) Markov inequality gives
$$
\sum_{i=1}^n \beta_{in}
 K_1^{(1)}
 \left(
\frac{\varepsilon_i-e}{b_1}
\right)
 =
 O_{\prob}
 \left( n
b_1^2 + (nb_1)^{1/2}
 \right)
 \max_{1 \leq i \leq n}
  \left|
\beta_{in}
 \right|,
$$
so that the proposition follows if we can prove that
\begin{equation}
\sup_{1 \leq i \leq n}
 \left|\beta_{in}\right|
 =
 O_{\prob}
\left( b_0^2\right),
 \label{BetasumTBP}
\end{equation}
as established now. For this, define
$$
\zeta_j (x)
 =
 \mathds{1}
 \left( x \in \mathcal{X}_0\right)
\left(m(X_j)-m(x)\right)
 K_0\left( \frac{X_j-x}{b_0}\right),
 \;\;
\nu_{in}(x)
 =
 \frac{1}{(n-1) b_0^d}
 \sum_{j=1, j\neq i}^n
 \left(
\zeta_j (x)-\esp[\zeta_j (x)]
 \right),
$$
and $\bar{\nu}_{n}(x)=\esp[\zeta_j(x)]/b_0^d $, so that
$$
\beta_{in}
=
 \frac{n-1}{n}\frac{\nu_{in} (X_i)
 +
 \bar{\nu}_n(X_i)}{\widehat{g}_{in}}
\;.
$$
For $\max_{1\leq i\leq n}|\bar{\nu}_n(X_i)|$, first observe
that  a second-order Taylor expansion applied successively to
$g(\cdot)$ and $m(\cdot)$ give, for $b_0$ small enough, and for
any $x$, $z$ in $\mathcal{X}$,
\begin{eqnarray*}
 \lefteqn
 {
 \left[ m(x+b_0z)-m(x)\right]
 g(x+b_0z)
 }
\\
& =&
\left[
 b_0 m^{(1)}(x) z
 +
 \frac{b_0^2}{2}
z m^{(2)}(x +\theta_1 b_0 z)z^{\top}
 \right]
 \left[
 g(x)
 +
 b_0 g^{(1)}(x) z
 +
 \frac{b_0^2}{2}
 z g^{(2)}(x +\theta_2 b_0z)z^{\top}
 \right],
  \end{eqnarray*}
for some $\theta_1 =\theta_1 (x,b_0 z)$ and $\theta_2 =\theta_2 (x,b_0
z)$ in $[0,1]$. Therefore, since $\int\!z K(z)dz=0$ under $(A_7)$,
it follows that, by $(A_1)$, $(A_2)$ and
 $(A_3)$,
\begin{eqnarray}
\nonumber
 \max_{1 \leq i \leq n}
|\bar{\nu}_n (X_i)|
&\leq&
\sup_{x\in\mathcal{X}_0}
|\bar{\nu}_n(x)|
 =
 \sup_{x \in
 \mathcal{X}_0}
 \left|
 \int
 \left(m ( x + b_0 z) - m(x)\right)
 K_0 (z) g(x+b_0z)
 dz
\right|
\\
&\leq&
 Cb_0^2.
 \label{Betasum1}
\end{eqnarray}
 Consider now the term  $\max_{1\leq i \leq n}|\nu_{in}(X_i)|$.
  The Bernstein inequality (see e.g. Serfling (2002)) and
$(A_4)$ give, for any $t>0$,
\begin{eqnarray*}
\prob
 \left(
 \max_{1 \leq i \leq n}
 | \nu_{in} (X_i)|
 \geq t
\right) &\leq &
 \sum_{i=1}^n
 \prob
 \left(
 | \nu_{in} (X_i) |
 \geq
t \right)
 \leq
 \sum_{i=1}^n
 \int
 \prob
 \left(
  | \nu_{in} (x) |
\geq t
 \left| X_i = x \right.
 \right) g (x)
  dx
\\
& \leq &
 2n \exp
 \left(
  - \frac{ (n-1) t^2 }
  { 2\sup_{x \in\mathcal{X}_0}
\Var (\zeta_j (x)/b_0^d) + \frac{4M}{3b_0^d} t}
 \right),
\end{eqnarray*}
where $M$ is such that $\sup_{x \in \mathcal{X}_0}|\zeta_j (x)|
\leq M$. The definition of $\mathcal{X}_0$ given in $(A_2)$,
$(A_3)$, $(A_7)$ and the standard Taylor expansion yield, for
$b_0$ small enough,
$$
\sup_{x \in \mathcal{X}_0}
 | \zeta_j (x) |
 \leq C b_0,
 \;\;\;
\sup_{x \in \mathcal{X}_0}
 \Var (\zeta_j (x)/b_0^d)
\leq \frac{1}{b_0^d}
 \sup_{x \in \mathcal{X}_0}
 \int
 \left( m(x +b_0 z) - m(x) \right)^2
 K_0^2 (z) g(x+b_0z) dz
  \leq
  \frac{C
b_0^2}{b_0^d},
$$
so that, for any $t \geq 0$,
$$
\prob
 \left(
 \max_{1 \leq i \leq n}
 | \nu_{in} (X_i) |
 \geq t
\right) \leq 2n
 \exp
 \left(
 - \frac{(n-1) b_0^d t^2 /b_0^2}{C + C t/b_0}
\right).
$$
This gives
$$
\prob
\left(
 \max_{1 \leq i \leq n}
 |\nu_{in} (X_i)|
 \geq
  \left(
\frac{b_0^2 \ln n}{ (n-1) b_0^d} \right)^{1/2} t
\right)
\leq 2n
\exp
\left(
 - \frac{ t^2 \ln n }
 {C + C t \left( \frac{\ln n}{
(n-1) b_0^d} \right)^{1/2} }
\right)
 =
  o(1),
$$
provided that $t$ is large enough and under $(A_9)$. It then
follows that
$$
\max_{1 \leq i \leq n} | \nu_{in} (X_i) |
 =
 O_{\prob}
 \left(
 \frac{b_0^2 \ln n}{ n b_0^d} \right)^{1/2}.
$$
This order, (\ref{Betasum1}) and Lemma \ref{Estig}  show that
(\ref{BetasumTBP}) is proved,  since $(b_0^2\ln
n/(nb_0^d))^{1/2}=O\left(b_0^2\right)$ under $(A_9)$, and that
$$
 \beta_{in}
  =
 \frac{n-1}{n} \frac{\nu_{in} (X_i)
 +
 \bar{\nu}_n (X_i)}{\widehat{g}_{in}}.
 \eop
$$

\subsection*{Proof of Proposition \ref{Sigsum}}

Assumption $(A_4)$ implies that  $\Sigma_{in}$ is independent of
$\varepsilon_i$, and that $\esp_n[\Sigma_{in}]=0$. This yields
\begin{eqnarray}
 \esp_n
 \left[
 \sum_{i=1}^n
 \Sigma_{in}
 K_1^{(1)}
 \left(
 \frac{\varepsilon_i-e}{b_1}
 \right)
 \right]
  = 0.
\label{Sigsum1}
\end{eqnarray}
Moreover,  observe that
\begin{eqnarray}
\nonumber
 \lefteqn{
 \Var_n
  \left[
 \sum_{i=1}^n
 \Sigma_{in}
 K_1^{(1)}
  \left(
\frac{\varepsilon_i-e}{b_1}
 \right)
 \right]
 }
\\
&=&
 \sum_{i=1}^n
 \Var_n
  \left[
 \Sigma_{in}
 K_1^{(1)}
 \left(
 \frac{\varepsilon_i-e}{b_1}
 \right)
 \right]
 +
 \sum_{i=1}^n
 \sum_{j=1\atop j\neq i}^n
 \Cov_n
 \left[
 \Sigma_{in}
 K_1^{(1)}
 \left(
 \frac{\varepsilon_{i}-e}{b_1}
 \right)
 ,
 \Sigma_{jn}
 K_1^{(1)}
 \left(
 \frac{\varepsilon_{j}-e}{b_1}
 \right)
 \right].
\label{VarSigm}
\end{eqnarray}
By Lemma \ref{MomderK}-(\ref{MomderK1}) and $(A_4)$, the first term above gives
\begin{eqnarray}
\nonumber \sum_{i=1}^n
 \Var_n
 \left[
 \Sigma_{in} K_1^{(1)}
 \left(
 \frac{\varepsilon_i-e}{b_1}
 \right)
 \right]
 &\leq&
 \sum_{i=1}^n
 \esp_n
\left[ \Sigma_{in}^2
 \right]
 \esp
 \left[
 K_1^{(1)}
  \left(
\frac{\varepsilon_i-e}{b_1} \right)^2
 \right]
\\\nonumber
&\leq&
 \frac{C b_1\sigma^2}{(nb_0^d)^2}
 \sum_{i=1}^n
 \sum_{j=1\atop j \neq i}^{n}
 \frac{\mathds{1}
 (X_i \in\mathcal{X}_0)}
{\widehat{g}_{in}^2}
 K_0^2
\left(\frac{X_j-X_i}{b_0}\right)
\\
 &\leq&
 \frac{C b_1\sigma^2}{nb_0^d}
 \sum_{i=1}^n
 \frac{\mathds{1}
 (X_i\in \mathcal{X}_0)
 \widetilde{g}_{in}}
{\widehat{g}_{in}^2},
 \label{VarSigm1}
\end{eqnarray}
where $\sigma^2=\Var(\varepsilon)$ and
$$
\widetilde{g}_{in}
 =
 \frac{1}{n b_0^d}
 \sum_{j=1,j \neq i}^{n}
 K_0^2\left( \frac{X_j-X_i}{b_0}\right).
$$
For the sum of  conditional covariances in (\ref{VarSigm}), write
\begin{eqnarray*}
 \lefteqn{
 \sum_{i=1}^n
 \sum_{j=1\atop j\neq i}^n
\Cov_n \left[
 \Sigma_{in}
 K_1^{(1)}
 \left(
 \frac{\varepsilon_{i}-e}{b_1}
 \right)
 ,
 \Sigma_{jn}
 K_1^{(1)}
 \left(
  \frac{\varepsilon_{j}-e}{b_1}
 \right)
\right]
 }
\\
&=&
 \sum_{i=1}^n
 \sum_{j=1\atop j\neq i}^n
\esp_n
 \left[
 \Sigma_{in}
 \Sigma_{jn}
 K_1^{(1)}
 \left(
\frac{\varepsilon_{i}-e}{b_1}
 \right)
 K_1^{(1)}
  \left(
\frac{\varepsilon_{j}-e}{b_1}
 \right)
  \right]
\\
& = &
 \sum_{i=1}^n
 \sum_{j=1\atop j\neq i}^n
 \frac{\mathds{1}(X_{i}\in\mathcal{X}_0)
\mathds{1}(X_{j}\in\mathcal{X}_0)}
 {(n b_0^d)^2 \widehat{g}_{in} \widehat{g}_{jn}}
 \sum_{k=1\atop k\neq i}^n
\sum_{\ell=1\atop \ell\neq j}^n
 K_0
 \left(
\frac{X_{k}-X_{i}}{b_0}
 \right)
 K_0
 \left(
\frac{X_{\ell}-X_{j}}{b_0} \right) \esp
 \left[
 \xi_{ki}
 \xi_{\ell j}
 \right],
\end{eqnarray*}
where $\xi_{ki}= \varepsilon_k K_1^{(1)}\left((\varepsilon_{i}-e)/b_1\right)$.
Moreover, under $(A_4)$, it is seen that for $k\neq\ell$,
$\esp[\xi_{ki}\xi_{\ell j}] =0 $ when $\Card\{i, j, k,
\ell\}\geq3$. Therefore
 the symmetry of $K_0$ implies that
\begin{eqnarray*}
 \lefteqn{
 \sum_{i=1}^n
 \sum_{j=1\atop j\neq i}^n
\Cov_n \left[ \Sigma_{in}
 K_1^{(1)}
  \left(
\frac{\varepsilon_{i}-e}{b_1}
 \right),
 \Sigma_{jn}
K_1^{(1)}
 \left(
 \frac{\varepsilon_{j}-e}{b_1}
 \right)
 \right]
 }
 &&
\\
& = &
\sum_{i=1}^n
 \sum_{j=1\atop j\neq i}^n
\frac{\mathds{1}(X_{i} \in \mathcal{X}_0) \mathds{1} (X_{j}
 \in
\mathcal{X}_0)}
 {(n b_0^d)^2 \widehat{g}_{in} \widehat{g}_{jn}}
K_0^2 \left( \frac{X_{j}-X_{i}}{b_0} \right)
 \esp^2
 \left[
\varepsilon K_1^{(1)}
\left(
\frac{\varepsilon-e}{b_1}
\right)
\right]
\\
&&
+
\sum_{i=1}^n
 \sum_{j=1\atop j\neq i}^n
 \frac{ \mathds{1} (X_{i}\in \mathcal{X}_0)
\mathds{1} (X_{j} \in \mathcal{X}_0)}
 {(nb_0^d)^2\widehat{g}_{in} \widehat{g}_{jn}}
 \sum_{k=1\atop k\neq  i, j}^n
 K_0 \left(\frac{X_{k}-X_{i}}{b_0}\right)
 K_0\left(\frac{X_{k}-X_{j}}{b_0}\right)
 \esp[\varepsilon^2]
\esp^2
 \left[ K_1^{(1)}
 \left(
 \frac{\varepsilon-e}{b_1}
\right)
\right].
\\
\end{eqnarray*}
Hence from Lemma \ref{Estig} and  Lemma \ref{MomderK}-(\ref{MomderK1}), we deduce
\begin{eqnarray*}
\nonumber
 \lefteqn{
  \left|
 \sum_{i=1}^n
 \sum_{j=1\atop j\neq i}^n
\Cov_n
 \left[
  \Sigma_{in}
 K_1^{(1)}
 \left(
\frac{\varepsilon_{i}-e}{b_1}
 \right)
 ,
  \Sigma_{jn}
 K_1^{(1)}
 \left(
 \frac{\varepsilon_{j}-e}{b_1}
\right)
\right]
\right|
 }
\\
  &=&
 O_{\prob}
 \left(
 \frac{b_1^4}{n b_0^d}
 \right)
 \sum_{i=1}^n
 \frac{\mathds{1}(X_i\in \mathcal{X}_0)\widetilde{g}_{in}}
 {\widehat{g}_{in}}
 +
 O_{\prob} (b_1^4)
 \sum_{i=1}^n
 \frac{\mathds{1}(X_i\in \mathcal{X}_0)g_{in}}
 {\widehat{g}_{in}},
\end{eqnarray*}
where  $\widetilde{g}_{in}$ is defined as  in (\ref{VarSigm})
and
$$
 g_{in}
 =
 \frac{1}{(n b_0^d)^2}
 \sum_{j=1\atop j\neq i}^n
 \sum_{k=1\atop k\neq j, i}^n
 K_0\left(\frac{X_k-X_i}{b_0}\right)
 K_0\left(\frac{X_k-X_j}{b_0}\right).
$$
Moreover, using Lemma \ref{Estig} and some technical details, it can be shown that
 \begin{eqnarray*}
  \sum_{i=1}^n
 \frac{\mathds{1}(X_i\in \mathcal{X}_0)g_{in}}
 {\widehat{g}_{in}}
 =
 O_{\prob}(1),
 \quad
 \sum_{i=1}^n
 \frac{\mathds{1}(X_i\in \mathcal{X}_0)\widetilde{g}_{in}}
 {\widehat{g}_{in}^k}
 =
 O_{\prob}(1),
 \;\; k=1,2.
 \end{eqnarray*}
Substituting these orders and
(\ref{VarSigm1}) in (\ref{VarSigm}), yields, for $b_1$ small enough,
\begin{eqnarray*}
 \Var_n
 \left[
 \sum_{i=1}^n
 \Sigma_{in}
 K_1^{(1)}
 \left(
 \frac{\varepsilon_i-e}{b_1}
 \right)
 \right]
 =
 O_{\prob}
\left( \frac{b_1}{b_0^d}
 +
 \frac{b_1^4}{b_0^d}
 +
 nb_1^4
 \right)
 =
 O_{\prob}
 \left(
 \frac{b_1}{b_0^d}
 +
 nb_1^4
 \right).
\end{eqnarray*}
Finally, this order, (\ref{Sigsum1}) and the Markov
inequality give
$$
\sum_{i=1}^n
 \Sigma_{in}
 K_1^{(1)}
 \left(
 \frac{\varepsilon_i-e}{b_1}
 \right)
 =
 O_{\prob}
 \left(
 \frac{b_1}{b_0^d}+nb_1^4
 \right)^{1/2}. \eop
$$

\subsection*{Proof of Proposition \ref{VarZn}}

Observe that Lemma \ref{Indep} yields that $\zeta_{in}$ and $\zeta_{jn}$ are
 independent given $X_1,\ldots,X_n$ for some $C>0$ such that
 $\|X_i-X_j\|\geq Cb_0$. Therefore
\begin{eqnarray}
\Var_n\left(\sum_{i=1}^n\zeta_{in}\right)
=
\sum_{i=1}^n\Var_n\left(\zeta_{in}\right)
+
\sum_{i=1}^n
\sum_{j=1\atop j\neq i}^n
\mathds{1}
\left(\left\|X_i - X_j\right\|<C b_0\right)
\Cov_n\left( \zeta_{in},\zeta_{jn}\right).
\label{VarTn}
\end{eqnarray}
 Let $\esp_{in}[\cdot]
=\esp_n[\cdot| X_1,\ldots,X_n,\varepsilon_k, k\neq i]$.
Since $\widehat{m}_{in} - m(X_i)$ depends only on
 $\left(X_1,\ldots,X_n,\varepsilon_k, k\neq i\right)$,
\begin{eqnarray*}
\sum_{i=1}^n
 \Var_n
 \left(\zeta_{in}\right)
  \leq
 \sum_{i=1}^n
 \esp_n
 \left[
 \zeta_{in}^2
 \right]
 =
 \sum_{i=1}^n
 \esp_n
 \left[
 \mathds{1}
 \left(X_i\in\mathcal{X}_0\right)
 \left(\widehat{m}_{in} - m(X_i)\right)^4
  \esp_{in}
  \left[
  K_1^{(2)}
  \left(\frac{\varepsilon_{i}-e}{b_1}\right)^2
  \right]
 \right],
 \end{eqnarray*}
with, using by  Lemma \ref{MomderK}-(\ref{MomderK1}),
\begin{eqnarray*}
 \esp_{in}
 \left[
 K_1^{(2)}
 \left(\frac{\varepsilon_{i}-e}{b_1}\right)^2
 \right]
 =
 \int
  K_1^{(2)}
 \left(
 \frac{\epsilon-e}{b_1}
 \right)^2
 f(\epsilon)
  d\epsilon
 \leq
 C b_1.
\end{eqnarray*}
 Therefore  Lemma \ref{BoundEspmchap} implies that
\begin{eqnarray*}
\sum_{i=1}^n
 \Var_n\left(\zeta_{in}\right)
 &\leq&
 Cb_1
 \sum_{i=1}^n
 \esp_n
 \left[
 \mathds{1}
 \left(X_i\in\mathcal{X}_0\right)
 (\widehat{m}_{in} - m(X_i))^4
 \right]
 \\
 &=&
 O_{\prob}\left(nb_1\right)
 \left(b_0^4+\frac{1}{nb_0^d}\right)^2.
\end{eqnarray*}
For the sum of the conditional covariances of (\ref{VarTn}), the order
is derived from the following
equalities:
\begin{eqnarray}
 \sum_{i=1}^n
 \sum_{j=1\atop j\neq i}^n
 \mathds{1}
 \left(
 \left\|X_i - X_j \right\|<C b_0
 \right)
 \esp_n\left[\zeta_{in}\right]
 \esp_n\left[\zeta_{jn}\right]
 &=&
 O_{\prob}
 \left(n^2b_0^db_1^6\right)
\left(b_0^4+ \frac{1}{nb_0^d}\right)^2,
  \label{Covzeta2}
  \\
  \sum_{i=1}^n
  \sum_{j=1\atop j\neq i}^n
  \mathds{1}
  \left(
  \left\|X_i - X_j\right\|<C b_0
  \right)
  \esp_n
  \left[
  \zeta_{in}
  \zeta_{jn}
  \right]
  &=&
  O_{\prob}
 \left(n^2b_0^db_1^{7/2}\right)
\left(b_0^4+ \frac{1}{nb_0^d}\right)^2.
 \label{Covzeta1}
\end{eqnarray}
 Indeed, since $b_1$ goes to $0$ under ${\rm (A_9)}$, the equalities above
 ensure that
\begin{eqnarray*}
\lefteqn{
 \sum_{i=1}^n
 \sum_{j=1\atop j\neq i}^n
 \mathds{1}
 \left(
 \left\|X_i - X_j \right\|<C b_0
 \right)
 \Cov_n
 \left(\zeta_{in},\zeta_{jn}\right)
 }
 \\
 &=&
O_{\prob}
 \left[
 \left(n^2b_0^db_1^6\right)
 \left(b_0^4+ \frac{1}{nb_0^d}\right)^2
 +
 \left(n^2b_0^db_1^{7/2}\right)
 \left(b_0^4+ \frac{1}{nb_0^d}\right)^2
\right]
\\
&=& O_{\prob}
 \left(n^2b_0^db_1^{7/2}\right)
 \left(b_0^4+ \frac{1}{nb_0^d}\right)^2.
\end{eqnarray*}
Combining this with the inequality above and  (\ref{VarTn}), and applying the
(conditional) Markov inequality, gives the desired result of the lemma.
 Hence,  it remains to prove (\ref{Covzeta2}) and
(\ref{Covzeta1}).  To this end, note that by  Lemma \ref{MomderK}-(\ref{MomderK2})
and the Cauchy-Schwartz inequality we have
 \begin{eqnarray*}
 \left|
 \esp_{n}
 \left[\zeta_{in}\right]
 \right|
 &=&
 \left|
 \esp_n
 \left[
 \mathds{1}
 \left(X_i\in\mathcal{X}_0\right)
 (\widehat{m}_{in} - m(X_i))^2
 \esp_{in}
 \left[
 K_1^{(2)}
 \left(
 \frac{\varepsilon_i-e}{b_1}
 \right)
 \right]
 \right]
 \right|
 \\
 &=&
 \left|
 \int
  K_1^{(2)}
 \left(
 \frac{\epsilon-e}{b_1}
 \right)
 f(\epsilon)
 d\epsilon
 \right|
 \times
 \esp_n
 \left[
 \mathds{1}
 \left(X_i\in\mathcal{X}_0\right)
 (\widehat{m}_{in} - m(X_i))^2
 \right]
 \\
 &\leq&
 Cb_1^3
 \biggl(
  \esp_n
 \left[
 \mathds{1}
 \left(X_i\in\mathcal{X}_0\right)
 (\widehat{m}_{in} - m(X_i))^4
 \right]
 \biggr)^{1/2},
 \end{eqnarray*}
uniformly in $i$, so that (by Lemma \ref{BoundEspmchap})
\begin{eqnarray*}
\sup_{1\leq i, j\leq n}
 \left|
 \esp_{n}
 \left[\zeta_{in}\right]
 \esp_n
\left[\zeta_{jn}\right]
 \right|
 &\leq&
 Cb_1^6
 \sup_{1\leq i \leq n}
 \esp_n
 \left[
 \mathds{1}
 \left(X_i\in\mathcal{X}_0\right)
 (\widehat{m}_{in} - m(X_i))^4
 \right]
 \\
 &=&
 O_{\prob}
 \left(b_1^6\right)
\left(b_0^4+ \frac{1}{nb_0^d}\right)^2.
\end{eqnarray*}
Therefore, since
$$
\sum_{i=1}^n
 \sum_{j=1\atop j\neq i}^n
 \mathds{1}
 \left(
 \| X_i - X_j \|< C b_0
 \right)
 =O_{\prob}(n^2b_0^d),
$$
this gives
\begin{eqnarray*}
  \sum_{i=1}^n
 \sum_{j=1\atop j\neq i}^n
 \mathds{1}
 \left(
 \| X_i - X_j \|< C b_0
 \right)
 \esp_{n}
 \left[\zeta_{in}\right]
 \esp_{n}
 \left[\zeta_{jn}\right]
  =
 O_{\prob}
 \left(n^2b_0^db_1^6\right)
\left(b_0^4+ \frac{1}{nb_0^d}\right)^2,
\end{eqnarray*}
which   proves (\ref{Covzeta2}).

 For (\ref{Covzeta1}), let $\beta_{in}$ and $\Sigma_{in}$  be  
 as in the statement of Proposition \ref{Betasum} and 
 Proposition \ref{Sigsum} respectively,
 and define
and
 $Z_{in}=\mathds{1}(X_i\in\mathcal{X}_0)(\widehat{m}_{in}-m(X_i))^2$.
This gives $Z_{in}=(\beta_{in}+\Sigma_{in})^2$, so that,
 for any $i\neq j$,
 \begin{eqnarray}
 \esp_n
 \left[\zeta_{in}\zeta_{jn}\right]
 =
 \esp_n
 \left[
 Z_{in}
  K_1^{(2)}
 \left(
 \frac{\varepsilon_j-e}{b_1}
 \right)
  \esp_{in}
 \left[
 Z_{jn}
  K_1^{(2)}
 \left(
 \frac{\varepsilon_i-e}{b_1}
 \right)
 \right]
 \right],
 \label{Prodzeta}
 \end{eqnarray}
where
\begin{eqnarray}
\nonumber
 \lefteqn{
 \esp_{in}
 \left[
 Z_{jn}
  K_1^{(2)}
 \left(
 \frac{\varepsilon_i-e}{b_1}
 \right)
 \right]
}
\\
&=&
 \beta_{jn}^2
 \esp_{in}
 \left[
 K_1^{(2)}
 \left(
 \frac{\varepsilon_i-e}{b_1}
 \right)
 \right]
 +
 2\beta_{jn}
 \esp_{in}
 \left[
 \Sigma_{jn}
 K_1^{(2)}
 \left(
 \frac{\varepsilon_i-e}{b_1}
 \right)
 \right]
 +
 \esp_{in}
 \left[
 \Sigma_{jn}^2
 K_1^{(2)}
 \left(
 \frac{\varepsilon_i-e}{b_1}
 \right)
 \right].
\label{Covzeta3}
\end{eqnarray}
 By Lemma \ref{MomderK}-(\ref{MomderK2}), the first term
above gives
\begin{eqnarray}
\left|
\beta_{jn}^2
 \esp_{in}
 \left[
 K_1^{(2)}
 \left(
 \frac{\varepsilon_i-e}{b_1}
 \right)
 \right]
 \right|
 \leq
 Cb_1^3
 \beta_{jn}^2.
 \label{Covzeta4}
\end{eqnarray}
Under  ${\rm (A_4)}$,  the second term of (\ref{Covzeta3}) equals
\begin{eqnarray*}
\lefteqn{
 \frac{2\beta_{jn}}
 {nb_0^d\widehat{g}_{jn}}
\sum_{k=1, k\neq j}^n
 K_0
 \left(
 \frac{X_k-X_j}{b_0}
 \right)
\esp_{in}
 \left[
 \varepsilon_k
   K_1^{(2)}
 \left(
 \frac{\varepsilon_i-e}{b_1}
 \right)
 \right]
 }
 \\
 &=&
 \frac{2\beta_{jn}}
 {nb_0^d\widehat{g}_{jn}}
 K_0
 \left(
 \frac{X_i-X_j}{b_0}
 \right)
\esp_{in}
 \left[
 \varepsilon_i
   K_1^{(2)}
 \left(
 \frac{\varepsilon_i-e}{b_1}
 \right)
 \right].
\end{eqnarray*}
Therefore, since $K_0$ is bounded under ${\rm (A_6)}$,  we have
(using Lemma \ref{MomderK}-(\ref{MomderK2}))
\begin{eqnarray}
\left|
2\beta_{jn}
\esp_{in}
 \left[
 \Sigma_{jn}
   K_1^{(2)}
 \left(
 \frac{\varepsilon_i-e}{b_1}
 \right)
 \right]
\right|
 \leq
\frac{Cb_1^3|\beta_{jn}|}
 {nb_0^d\widehat{g}_{jn}}.
 \label{Covzeta5}
\end{eqnarray}
For the  last term of (\ref{Covzeta3}), write
\begin{eqnarray*}
\lefteqn{
 \esp_{in}
 \left[
 \Sigma_{jn}^2
   K_1^{(2)}
 \left(\frac{\varepsilon_i-e}{b_1}\right)
 \right]
}
\\
&=&
\frac{\mathds{1}(X_j\in\mathcal{X}_0)}{(nb_0^d\widehat{g}_{jn})^2}
 \sum_{k=1\atop k\neq j}^n
\sum_{\ell=1\atop\ell\neq j}^n
 K_0\left(\frac{X_k-X_j}{b_0}\right)
 K_0\left(\frac{X_{\ell}-X_j}{b_0}\right)
 \esp_{in}
 \left[
 \varepsilon_k
 \varepsilon_{\ell}
  K_1^{(2)}
 \left(\frac{\varepsilon_i-e}{b_1}\right)
 \right]
 \\
 &=&
 \frac{\mathds{1}(X_j\in\mathcal{X}_0)}{(nb_0^d\widehat{g}_{jn})^2}
 \sum_{k=1, k\neq j}^n
 K_0^2\left(\frac{X_k-X_j}{b_0}\right)
 \esp_{in}
 \left[
 \varepsilon_k^2
  K_1^{(2)}
 \left(\frac{\varepsilon_i-e}{b_1}\right)
 \right].
\end{eqnarray*}
Since
\begin{eqnarray*}
\lefteqn
 {
 \left|
 \esp_{in}
 \left[
 \varepsilon_k^2
  K_1^{(2)}
 \left(\frac{\varepsilon_i-e}{b_1}\right)
 \right]
 \right|
}
\\
&\leq& \max
 \left\lbrace
\sup_{e^{\prime}\in\Rit}
\left|
 \esp_{in}
 \left[
 \varepsilon^2
  K_1^{(2)}
 \left(\frac{\varepsilon-e^{\prime}}{b_1}\right)
 \right]
 \right|,
 \;
 \esp[\varepsilon^2]
\sup_{e^{\prime}\in\Rit}
 \left|
 \esp_{in}
 \left[
 K_1^{(2)}
 \left(\frac{\varepsilon-e^{\prime}}{b_1}\right)
 \right]
 \right|
\right\rbrace
\\
&\leq&
 C b_1^3,
\end{eqnarray*}
uniformly in $i$, this gives
$$
 \left|
 \esp_{in}
 \left[
 \Sigma_{jn}^2
   K_1^{(2)}
 \left(
 \frac{\varepsilon_i-e}{b_1}
 \right)
 \right]
 \right|
 \leq
\frac{Cb_1^3\mathds{1}(X_j\in\mathcal{X}_0)}{(nb_0^d\widehat{g}_{jn})^2}
\sum_{k=1, k\neq j}^n
K_0^2\left(\frac{X_k-X_j}{b_0}\right).
$$
Substituting this,  (\ref{Covzeta5}) and  (\ref{Covzeta4}) in
(\ref{Covzeta3}), yields
$$
\left|
 \esp_{in}
 \left[
 Z_{jn}
  K_1^{(2)}
 \left(
 \frac{\varepsilon_i-e}{b_1}
 \right)
 \right]
\right|
 \leq
 Cb_1^3 M_n,
$$
where
$$
M_n
=
\sup_{1\leq j\leq n}
\left[
\beta_{jn}^2
+
\frac{|\beta_{jn}|}{nb_0^d\widehat{g}_{jn}}
+
\frac{\mathds{1}(X_j\in\mathcal{X}_0)}{(nb_0^d\widehat{g}_{jn})^2}
\sum_{k=1, k\neq j}^n
K_0^2\left(\frac{X_k-X_j}{b_0}\right)
 \right].
$$
Hence from (\ref{Prodzeta}), the Cauchy-Schwarz inequality
and Lemmas \ref{BoundEspmchap}, \ref{MomderK} we deduce
\begin{eqnarray*}
\lefteqn{
\sum_{i=1}^n
\sum_{j=1\atop j\neq i}^n
 \mathds{1}
 \left(
 \| X_i - X_j \|< C b_0
 \right)
\left|
 \esp_n
 \left[\zeta_{in}\zeta_{jn}\right]
 \right|
 }
\\
&\leq&
Cb_1^3M_n
\sum_{i=1}^n
\sum_{j=1\atop j\neq i}^n
 \mathds{1}
 \left(
 \| X_i - X_j \|< C b_0
 \right)
\esp_n
 \left|
 Z_{in} K_1^{(2)}
\left(\frac{\varepsilon_j-e}{b_1}\right)
\right|
\\
&\leq&
Cb_1^3M_n
 \sum_{i=1}^n
 \sum_{j=1\atop j\neq i}^n
 \mathds{1}
 \left(
 \| X_i - X_j \|< C b_0
 \right)
\esp_n^{1/2}
\left[
\mathds{1}(X_i\in\mathcal{X}_0)
(\widehat{m}_{in}-m(X_i))^4
\right]
 \esp_n^{1/2}
 \left[
 K_1^{(2)}
 \left(
 \frac{\varepsilon_j-e}{b_1}
 \right)^2
 \right]
 \\
 &=&
b_1^{7/2} M_n
 O_{\prob}
 \left(b_0^4+\frac{1}{nb_0^d}\right)
\sum_{i=1}^n \sum_{j=1\atop j\neq i}^n
\mathds{1}
\left(\|X_i-X_j\|\leq Cb_0\right).
\end{eqnarray*}
Moreover, using Proposition \ref{Betasum} (which gives 
$\sup_{i}|\beta_{in}|=O_{\prob}(b_0^2)$),
Lemma \ref{Estig} and some technical details,
it can be shown that
$$
M_n =
 O_{\prob}
 \left(
 b_0^4
 +
  \frac{b_0^2}{nb_0^d}
  +
\frac{1}{nb_0^d}
 \right)
= O_{\prob}
\left(b_0^4+ \frac{1}{nb_0^d}\right).
$$
Substituting this order in the inequality above, yields
 (\ref{Covzeta1}) and  finishes  the proof of the
proposition. \eop

\subsection*{Proof of Proposition \ref{VarRn}}

Observe that by Lemma \ref{Indep}, we have
\begin{eqnarray}
\Var_n\left(\sum_{i=1}^n R_{in}\right)
=
\sum_{i=1}^n\Var_n\left(R_{in}\right)
+
\sum_{i=1}^n
\sum_{j=1\atop j\neq i}^n
\mathds{1}
\left(\left\|X_i - X_j\right\|<C b_0\right)
\Cov_n
\left( R_{in},R_{jn}
 \right).
 \label{VarRn1}
\end{eqnarray}
Let $\esp_{in}[\cdot]
=\esp_n[\cdot| X_1,\ldots,X_n,\varepsilon_k, k\neq i]$, and write
\begin{eqnarray*}
\sum_{i=1}^n
 \Var_n
 \left(R_{in}\right)
  \leq
 \sum_{i=1}^n
 \esp_n
 \left[
 R_{in}^2
 \right]
 =
 \sum_{i=1}^n
 \esp_{n}
 \biggl[
 \mathds{1}
 \left( X_i \in \mathcal{X}_0\right)
 (\widehat{m}_{in} - m(X_i))^6
 \esp_{in}
 \left[ I_{in}^2\right]
 \biggr],
\end{eqnarray*}
with, using  ${\rm (A_4)}$, the Cauchy-Schwarz inequality and
Lemma \ref{MomderK}-(\ref{MomderK3}),
\begin{eqnarray*}
\lefteqn{
 \esp_{in}
 \left[ I_{in}^2\right]
 =
\esp_{in}
 \left[
\left\lbrace
\int_{0}^{1}
 (1-t)^2
 K_1^{(3)}
 \left(
 \frac{
 \varepsilon_i-t(\widehat{m}_{in}-m(X_i))-e}{b_1}
 \right)
 dt
 \right\rbrace^2
\right]
 }
\\
 &\leq&
 \esp_{in}
 \left[
\int_{0}^{1}
(1-t)^4
 K_1^{(3)}
 \left(
 \frac{\varepsilon_i-t(\widehat{m}_{in} -m(X_i))-e}{b_1}
 \right)^2
 dt
 \right]
\\
&=&
 \int_{0}^{1}
 (1-t)^4
 \left[
 \int
 K_1^{(3)}
 \left(
 \frac{\epsilon-t(\widehat{m}_{in} -m(X_i))-e}{b_1}
 \right)^2
 f(\epsilon)
 d\epsilon
 \right]
 dt
\leq
 Cb_1.
\end{eqnarray*}
Therefore Lemma \ref{BoundEspmchap} implies that
\begin{eqnarray}
\nonumber
 \sum_{i=1}^n
 \Var_n\left(R_{in}\right)
 &\leq&
 C nb_1
 \sup_{1\leq i\leq n}
 \esp_{n}
 \left[
 \mathds{1}
 \left( X_i \in \mathcal{X}_0\right)
 (\widehat{m}_{in} - m(X_i))^6
 \right]
 \\
&=&
 O_{\prob}
 \left(nb_1\right)
 \left(b_0^4+\frac{1}{nb_0^d}\right)^3.
 \label{VarRn2}
\end{eqnarray}
 For the second  term of (\ref{VarRn1}), write
\begin{eqnarray*}
\lefteqn{
\left|
\Cov_n\left(R_{in}, R_{jn}\right)
 \right|
 \leq
 \left(
 \Var_n\left(R_{in}\right)
 \Var_n\left(R_{jn}\right)
 \right)^{1/2}
 }
 \\
 &\leq&
 Cb_1
 \sup_{1\leq i\leq n}
 \esp_{n}
 \left[
 \mathds{1}
 \left( X_i \in \mathcal{X}_0\right)
 (\widehat{m}_{in} - m(X_i))^6
 \right]
 =
  O_{\prob}(b_1)
 \left(b_0^4+\frac{1}{nb_0^d}\right)^3,
\end{eqnarray*}
uniformly in $i$ and $j$. Hence
\begin{eqnarray*}
\lefteqn{
\sum_{i=1}^n
\sum_{j=1\atop j\neq i}^n
\left(\|X_i-X_j\|\leq Cb_0\right)
\left|
\Cov_n\left(R_{in}, R_{jn}\right)
\right|
}
\\
&=&
 O_{\prob}\left(b_1\right)
 \left(b_0^4+\frac{1}{nb_0^d}\right)^3
 \sum_{i=1}^n
 \sum_{j=1\atop j\neq i}^n
 \left(\|X_i-X_j\|\leq Cb_0\right)
 \\
 &=&
 O_{\prob}\left(b_1\right)
 \left(b_0^4+\frac{1}{nb_0^d}\right)^3
 \left(n^2b_0^d\right).
\end{eqnarray*}
Finally, this order, (\ref{VarRn2}) and (\ref{VarRn1}) give, since $nb_0^d$
diverges  under ${\rm (A_8)}$,
$$
\Var\left(\sum_{i=1}^nR_{in}\right)
=
O_{\prob}\left(nb_1+n^2b_0^db_1\right)
 \left(b_0^4+\frac{1}{nb_0^d}\right)^3
 =
 O_{\prob}
  \left(n^2b_0^db_1\right)
 \left(b_0^4+\frac{1}{nb_0^d}\right)^3.
 \eop
$$

\setcounter{equation}{0}
\setcounter{subsection}{0}
\renewcommand{\theequation}{\arabic{equation}}

\begin{center}
\section*{Appendix:  Proof of the  intermediate results}
\end{center}

\subsection*{Proof of Lemma \ref{Estig}}

First note that  by $(A_7)$, we have
 $
 \int\!
z K_0(z) dz
 =0
 $
 and
 $
 \int\!
K_0(z) dz
 =1$.  Therefore
$(A_1)$, $(A_2)$ and the second-order Taylor expansion,
 yield, for  $b_0$ small enough and any $x$ in $\mathcal{X}_0$,
\begin{eqnarray*}
 \lefteqn{
 \left|
 \overline{g}_{n}(x)-g(x)
 \right|
  =
 \left|
 \frac{1}{b_0^d}
 \int
 K_0
 \left(\frac{z-x}{b_0}\right)
 g(z)
 dz
 -
 g(x)
 \right|
  =
 \left|
\int K_0(z) \left[g(x+b_0z)-g(x)\right]
 dz
 \right|
 }
 &&
\\
&=&
 \left|
 \int K_0(z)
 \left[
 b_0 g^{(1)}(x) z
 +
 \frac{b_0^2}{2}
 z g^{(2)} (x + \theta b_0 z)z^{\top}
\right] dz
 \right|,
 \;\theta = \theta (x,b_0 z)\in [0,1]
\\
& = & \frac{b_0^2}{2}
 \left|
 \int z g^{(2)}(x+\theta b_0z)
 z^{\top} K_0(z) dz
 \right|
 \leq C b_0^2,
\end{eqnarray*}
so that
$$
\sup_{x\in\mathcal{X}_0}
\left|\overline{g}_n (x)-g(x)\right|
=
O\left(b_0^2\right).
$$
This gives the first result of the lemma. To prove the second
and third results of  the lemma, note that it is
 sufficient to show that
$$
 \sup_{x\in\mathcal{X}_0}
 \left|
 \widehat{g}_{n}(x)
  -
 \overline{g}_n(x)
 \right|
 =
 O_{\prob}
 \left(
  \frac{\ln n}{nb_0^d}
 \right)^{1/2},
$$
 since $\bar{g}_n (x)$ is asymptotically bounded away from
$0$ over
 $\mathcal{X}_0$
 and that
$|\overline{g}_n (x) - g(x) | = O (b_0^2) $ uniformly for $x$ in
$\mathcal{X}_0$. This follows from Theorem 1 in  Einmahl and Mason
(2005). \eop

\subsection*{Proof of Lemma \ref{Indep}}

Since $K_0(\cdot)$ has a compact support under ${\rm (A_6)}$,
there is a $C>0$ such that $\| X_i - X_j \| \geq C b_0$ implies
that for any integer number $k$ of $[1,n]$, $K_0 ( (X_k -
X_i)/b_0) = 0$ if $K_0 ( (X_j - X_k)/b_0) \neq 0$. Let $D_j
\subset [1,n]$ be such that an integer number $k$ of $[1,n]$ is in
$D_j$ if and only if $K_0 ( (X_j - X_k)/b_0) \neq 0$. Abbreviate
$\prob (\cdot| X_1, \ldots,X_n)$ into $\prob_n$ and assume that
$\| X_i - X_j \| \geq C b_0$ so that $D_i$ and $D_j$ have an empty
intersection. Note also that taking $C$ large enough ensures that
$i$ is not in $D_j$ and $j$ is not in $D_i$. It then follows,
under ${\rm (A_4)}$ and since $D_i$ and $D_j$ only depend upon
$X_1,\ldots,X_n$,
\begin{eqnarray*}
\lefteqn{
 \prob_n
 \biggl(
 \left(
 \widehat{m}_{in} - m(X_i),\varepsilon_i
 \right)
 \in A \mbox{ \rm and }
 \left(
\widehat{m}_{jn} - m (X_j),\varepsilon_j
\right)
 \in B \biggr)
 }
&&
\\
& = &
\prob_n \left( \left( \frac{\sum_{k \in D_i \setminus\{i\}}
 \left( m(X_{k}) - m(X_i) + \varepsilon_{k} \right)
  K_0
\left( (X_{k} - X_i)/b_0 \right)} {\sum_{k \in D_i \setminus\{i\}}
 K_0 \left( (X_{k} - X_i)/b_0 \right)}
 ,
 \varepsilon_i\right) \in A \right.
\\
&& \;\;\;\;\;\;\;\;\;\;\;\;\;\;\;\;\;\;\;\; \left. \mbox{ \rm and}
\left( \frac{ \sum_{\ell \in D_j \setminus \{j \}}
 \left(
m(X_{\ell}) - m(X_j) + \varepsilon_{\ell} \right)
 K_0 \left(
(X_{\ell} - X_j)/b_0 \right) } { \sum_{\ell \in D_j \setminus \{j
\}} K_0 \left( (X_{\ell} - X_j)/b_0 \right)} ,
\varepsilon_j\right) \in B \right)
\\
& = & \prob_n \left( \left( \frac{ \sum_{k \in D_i \setminus \{i
\}} \left( m(X_{k}) - m(X_i) + \varepsilon_{k}\right)
 K_0 \left(
(X_{k} - X_i)/b_0 \right)} {\sum_{k \in D_i \setminus \{i \}}
K_0\left( (X_{k} - X_i)/b_0 \right)} ,
 \varepsilon_i \right) \in
A \right)
\\
&& \;\;\;\;\;\;\;\;\;\;\;\;\;\;\;\;\;\;\;\;
 \times\;
 \prob_n
 \left(
\left( \frac{ \sum_{\ell \in D_j \setminus \{j \}} \left(
m(X_{\ell}) - m(X_j) + \varepsilon_{\ell} \right) K_0 \left(
(X_{\ell} - X_j)/b_0 \right) } { \sum_{\ell \in D_j \setminus \{j
\}} K_0\left( (X_{\ell} - X_j)/b_0\right) } , \varepsilon_j
\right) \in B \right)
\\
& = & \prob_n
 \left(
 \left(\widehat{m}_{in} - m(X_i), \varepsilon_i \right)
 \in A
 \right)
\times \prob_n \left( \left(\widehat{m}_{jn} - m
(X_j),\varepsilon_j \right)
 \in B
\right).
\end{eqnarray*}
This gives the result of Lemma \ref{Indep}, since both
$\left(\widehat{m}_{in} - m (X_i), \varepsilon_i\right)$ and
$\left(\widehat{m}_{jn} - m (X_j), \varepsilon_j\right)$ are
independent given $X_1, \ldots, X_n$. \eop

\subsection*{Proof of Lemma \ref{BoundEspmchap}}

Let $\beta_{in}$  as in the statement of Proposition \ref{Betasum} and set
\begin{eqnarray*}
g_{in}
 =
 \frac{1}{nb_0^d}
 \sum_{j=1, j\neq i}^n
 K_0^4\left(\frac{X_j-X_i}{b_0}\right),
 \quad
\widetilde{g}_{in}
 =
 \frac{1}{nb_0^d}
 \sum_{j=1, j\neq i}^n
 K_0^2\left(\frac{X_j-X_i}{b_0}\right).
\end{eqnarray*}
The proof of the lemma  is based on the following bound:
\begin{eqnarray}
\esp_n
 \left[
 \mathds{1}
 \left(X_i\in\mathcal{X}_0\right)
 \left(\widehat{m}_{in} - m(X_i)\right)^k
\right]
 \leq
 C
 \left[
 \beta_{in}^k
 +
 \frac{
 \mathds{1}
 \left(X_i\in\mathcal{X}_0\right)
 \widetilde{g}_{in}^{k/2}}
 {(nb_0^d)^{(k/2)}\widehat{g}_{in}^k}
 \right],
 \quad
 k\in\{4,6\}.
 \label{Espm}
 \end{eqnarray}
 Indeed, taking successively $k=4$ and $k=6$ in (\ref{Espm}),
 we have, by (\ref{BetasumTBP}) and Lemma \ref{Estig}
\begin{eqnarray*}
 \sup_{1\leq i\leq n}
 \esp_n
 \left[
 \mathds{1}\left(X_i\in\mathcal{X}_0\right)
 \left(\widehat{m}_{in} - m(X_i)\right)^4
\right]
 &=&
 O_{\prob}
 \left(
 b_0^8
  +
 \frac{1}{(nb_0^d)^2}
 \right)
 =
 O_{\prob}
 \left(b_0^4+\frac{1}{nb_0^d}\right)^2,
 \\
 \sup_{1\leq i\leq n}
 \esp_n
 \left[
 \mathds{1}
 \left(X_i\in\mathcal{X}_0\right)
 \left(\widehat{m}_{in} - m(X_i)\right)^6
\right]
 &=&
 O_{\prob}
 \left(
 b_0^{12}
  +
 \frac{1}{(nb_0^d)^3}
 \right)
 =
 O_{\prob}
 \left(b_0^4+\frac{1}{nb_0^d}\right)^3,
\end{eqnarray*}
which gives the desired results.
 Hence it remains to prove (\ref{Espm}).
 To this end,  let $\Sigma_{in}$ be as in the statement of Proposition \ref{Sigsum}
 and observe that
 $\mathds{1}(X_i \in\mathcal{X}_0)
 \left(\widehat{m}_{in}- m(X_i)\right)
 =\beta_{in}+\Sigma_{in}$. Since
 $\beta_{in}$ depends only on
 $\left(X_1,\ldots,X_n\right)$, this gives, for any $k\in\{4,6\}$,
 we have
\begin{eqnarray}
\esp_n
 \left[
 \mathds{1}
 (X_i \in \mathcal{X}_0)
 \left(\widehat{m}_{in}- m(X_i)\right)^k
 \right]
  \leq
 C\beta_{in}^k
 +
 C\esp_n\left[\Sigma_{in}^k\right].
 \label{Espm5}
\end{eqnarray}
The order  of the second term of  (\ref{Espm5}) is computed
by applying Theorem 2 in Whittle (1960) or the
Marcinkiewicz-Zygmund inequality (see e.g Chow and Teicher, 2003,
p. 386). These inequalities show that for linear form
$L=\sum_{j=1}^n a_j\zeta_j$ with independent mean-zero random
variables
 $\zeta_1,\ldots,\zeta_n$, it holds that, for any $k\geq 1$,
 $$
 \esp
 \left|L^k\right|
 \leq
  C(k)
  \left[
 \sum_{j=1}^n
  a_j^2
\esp^{2/k}
\left|\zeta_j^k \right|
 \right]^{k/2},
 $$
where $C(k)$ is a positive real depending only on $k$. Now,
observe that for any integer $i\in[1,n]$,
$$
\Sigma_{in}
 =
 \sum_{j=1, j\neq i}^n
 \sigma_{jin},
 \quad
 \sigma_{jin}
 =
 \frac{\mathds{1}
 \left( X_i \in \mathcal{X}_0\right)}
 {nb_0^d\widehat{g}_{in}}
 \varepsilon_j
 K_0\left(\frac{X_j-X_i}{b_0}\right).
$$
Since under $(A_4)$, the $\sigma_{jin}$'s, $(j=1,\ldots,n)$, are
centered independent variables given $X_1,\ldots,X_n$,  this
yields, for any $k\in\{4,6\}$,
\begin{eqnarray*}
 \esp_n
 \left[
 \Sigma_{in}^k\right]
 \leq
 C\esp\left[\varepsilon^k\right]
 \left[
 \frac{\mathds{1}
 \left( X_i \in \mathcal{X}_0\right)}
 {(nb_0^d)^2\widehat{g}_{in}^2}
 \sum_{j=1}^n
  K_0^2
 \left(
 \frac{X_j-X_i}{b_0}
 \right)
 \right]^{k/2}
 \leq
 \frac{C\mathds{1}
 \left(X_i\in\mathcal{X}_0\right)\widetilde{g}_{in}^{k/2}}
 {(nb_0^d)^{(k/2)}\widehat{g}_{in}^k}\;.
\end{eqnarray*}
 Hence, this  bound and (\ref{Espm5})  imply that
$$
 \esp_n
 \left[
 \mathds{1}(X_i \in \mathcal{X}_0)
 \left(\widehat{m}_{in}- m(X_i)\right)^k
 \right]
 \leq C\left[
 \beta_{in}^k
  +
 \frac{
 \mathds{1}
 \left(X_i\in\mathcal{X}_0\right)
 \widetilde{g}_{in}^{k/2}}
 {(nb_0^d)^{(k/2)}\widehat{g}_{in}^k}
 \right],
$$
which proves (\ref{Espm}), and then completes the proof of the
lemma. \eop

\subsection*{Proof of Lemma \ref{MomderK}}

Set $h_p(e)=e^p f(e)$, $p\in[0,2]$.
For  the first inequality of
(\ref{MomderK1}), note that  under ${\rm (A_5)}$ and ${\rm
(A_7)}$, the change of variable $\epsilon=e+b_1 v$ give, for any
$\ell\in\{1,2,3\}$,
\begin{eqnarray}
\nonumber
 \left|
 \int
 K_1^{(\ell)}
 \left(\frac{\epsilon-e}{b_1}\right)^2
 \epsilon^p f(\epsilon)
 d\epsilon
 \right|
 &=&
 \left|
 b_1
 \int
 K_1^{(\ell)}(v)^2
 h_p(e+b_1v)
 dv
 \right|
 \\\nonumber
 &\leq&
 b_1
 \sup_{e\in\Rit}
 |h_p(e)|
 \int
 K_1^{(\ell)}(v)^2
 dv
 \\
 &\leq&
 Cb_1,
 \label{Ineg1}
\end{eqnarray}
which yields the first inequality of (\ref{MomderK1}). For the
second inequality of (\ref{MomderK1}), observe that under ${\rm
(A_7)}$ we have $\int \! K_1^{(\ell)}(v)dv =0$. Therefore, since
$h_p(\cdot)$ has bounded second order derivatives under ${\rm
(A_5)}$, the Taylor inequality gives
\begin{eqnarray*}
\left| \int
 K_1^{(\ell)}
\left(\frac{\epsilon-e}{b_1}\right)
 \epsilon^p f(\epsilon)
 d\epsilon
 \right|
 &=&
  b_1
 \left|
 \int
 K_1^{(\ell)}(v)
 \left[
 h_p(e+b_1v)-h_p(e)
 \right]
 \right|
  dv
 \\
 &\leq&
 b_1^2
 \sup_{e\in\Rit}|h_p^{(1)}(e)|
 \int
 |v K_1^{(\ell)}(v)|
  dv
 \leq
 Cb_1^2.
\end{eqnarray*}
which completes the  proof of (\ref{MomderK1}).
 The first inequalities of  (\ref{MomderK2}) and
 (\ref{MomderK3}) are given by (\ref{Ineg1}). The second inequalities of
(\ref{MomderK2}) and (\ref{MomderK3}) are proved simultaneously.
For this, note that for any integer $\ell\in\{2,3\}$,
$$
\int
 K_1^{(\ell)}
 \left(\frac{\epsilon-e}{b_1}\right)
 h_p(\epsilon)
 d\epsilon
 =
b_1
\int
 K_1^{(\ell)}(v)
h_p(e+b_1v) dv.
$$
By ${\rm (A_7)}$,  $K_1$ has  a compact support and  satisfies
$\int \! K_1^{(\ell)}(v)dv=0$ and
 $\int\! v K_1^{(\ell)}(v) dv=0$.
 Hence the second order Taylor expansion applied
 to $h_p(\cdot)$ gives, for some
$\theta=\theta (e,b_1, v)\in[0,1]$,
\begin{eqnarray*}
 \lefteqn{
 \left|
 \int
 K_1^{(\ell)}
\left(
 \frac{\epsilon-e}{b_1}
 \right)
  h_p(\epsilon)
  d\epsilon
 \right|
 =
 \left|
  b_1
 \int
 K_1^{(\ell)}(v)
\left[
 h_p(e+b_1v)
 -
  h_p(e)
 \right]
 dv
 \right|
 }
 &&
\\
&=&
\left|
b_1
\int
 K_1^{(\ell)}(v)
\left[
 b_1 v h_p^{(1)}(e)
 +
  \frac{b_1^2v^2}{2}
h_p^{(2)}(e+\theta b_1v) \right]
 dv
 \right|
 \\
  &=&
 \left|
\frac{b_1^3}{2}
\int
v^2 K_1^{(\ell)}(v)
 h_p^{(2)}(e+\theta b_1v)
 dv
\right|
\\
&\leq&
 \frac{b_1^3}{2}
 \sup_{e\in \Rit}|h_p^{(2)}(e)|
 \int
\left| v^2K_1^{(\ell)}(v)
 \right|
 dv
 \leq
 Cb_1^3.\eop
\end{eqnarray*}

\end{document}